\documentclass[12pt]{amsart}
\usepackage{amsmath,amssymb,amsthm,graphics,amscd,mathrsfs}
\usepackage{amscd, amsfonts, amssymb, amsthm}
\usepackage{fullpage, verbatim}
\usepackage{hyperref} 
\usepackage{breakurl}
\usepackage{array}
\usepackage{latexsym}
\usepackage{enumerate}

\newcommand\CA{{\mathcal A}}

\newcommand\CIF{{\mathcal {IF}}}

\newcommand\CI{{\mathcal I}}

\newcommand\BBR{{\mathbb R}}
\newcommand\BBZ{{\mathbb Z}}

\newcommand\codim{\operatorname{codim}}

\newcommand\Der{{\operatorname{Der}}}

\newcommand\GL{\operatorname{GL}}

\newcommand\pdeg{\operatorname{pdeg}}

\newcommand\hgt{\operatorname{ht}}

\DeclareMathOperator{\Aut}{Aut}

\numberwithin{equation}{section}

\theoremstyle{plain}

\newtheorem{theorem}[equation]{Theorem}

\theoremstyle{definition}
\newtheorem{defn}[equation]{Definition}
\newtheorem{remark}[equation]{Remark}
\newtheorem{remarks}[equation]{Remarks}
\newtheorem{example}[equation]{Example}

\newtheorem{algor}[equation]{Algorithm}

\newcommand{\algo}[6]
{
\begin{algor}{{\tt #1}{\rm (#2)}}\label{#1}\end{algor}
\vspace{-6pt}\noindent{\it #3}.

{\bf Input:} #4

{\bf Output:} #5

\newcounter{#1}
\begin{list}{\textbf{\arabic{#1}.}}{\usecounter{#1}}
#6\end{list}\vspace{3pt}}

\newcommand{\RS}{R}
\newcommand{\SI}{\Delta}

\subjclass[2010]{20F55, 52B30, 52C35, 14N20}

\begin{document}

\title[Arrangements of ideal type are inductively free]
{Arrangements of ideal type\\ are inductively free}

\author[M. Cuntz]{Michael Cuntz}
\address
{Institut f\"ur Algebra,~Zahlentheorie und Diskrete Mathematik,
Fakult\"at f\"ur Mathematik und Physik,
Gottfried Wilhelm Leibniz Universit\"at Hannover,
Welfengarten 1, D-30167 Hannover, Germany}
\email{cuntz@math.uni-hannover.de}

\author[G. R\"ohrle]{Gerhard R\"ohrle}
\address
{Fakult\"at f\"ur Mathematik,
Ruhr-Universit\"at Bochum,
D-44780 Bochum, Germany}
\email{gerhard.roehrle@rub.de}

\author[A. Schauenburg]{Anne Schauenburg}
\address
{Fakult\"at f\"ur Mathematik,
Ruhr-Universit\"at Bochum,
D-44780 Bochum, Germany}
\email{anne.schauenburg@rub.de}


\keywords{Root system, Weyl arrangement, Weyl, groupoid, arrangement of ideal type, free arrangement, inductively free arrangement}

\allowdisplaybreaks

\begin{abstract}
Extending earlier work by Sommers and Tymoczko, in 2016 Abe, Barakat, Cuntz, Hoge, and Terao established that each arrangement of ideal type $\mathcal{A}_\mathcal{I}$ stemming from an ideal $\mathcal{I}$ in the set of positive roots of a reduced root system is free.

Recently, R\"ohrle showed that a large class of the $\mathcal{A}_\mathcal{I}$ satisfy the stronger property of inductive freeness and conjectured that this property holds for all $\mathcal{A}_\mathcal{I}$. In this article, we confirm this conjecture.
\end{abstract}

\maketitle

\section{Introduction}

In this paper we study certain arrangements which are 
associated with ideals in the set of positive roots of 
a reduced root system, so called 
\emph{arrangements of ideal type} $\CA_\CI$,
Definition \ref{def:idealtype},
cf.~\cite[\S 11]{ST06}.
Our goal is to confirm a conjecture from 
\cite[Conj.~1.16]{gR17}
that all $\CA_\CI$ are inductively free, 
Definition \ref{def:indfree}.

\subsection{Ideals in $\RS_+$}
\label{sect:ideals}
Let $\RS$ be an irreducible, reduced root system
and let $\RS_+$ be the set of positive roots 
with respect to some set of simple roots $\SI$.
An \emph{(upper) order ideal}, 
or simply \emph{ideal} for short, 
of  $\RS_+$, is a subset $\CI$  of $\RS_+$ 
satisfying the following condition: 
if $\alpha \in \CI$ and $\beta \in \RS_+$ so that 
$\alpha + \beta \in \RS_+$, then $\alpha + \beta \in \CI$.

Recall the standard partial ordering 
$\preceq$ on $\RS$: $\alpha \preceq \beta$
provided $\beta - \alpha$ is a $\BBZ_{\ge0}$-linear combination 
of positive roots, or $\beta = \alpha$. Then $\CI$ is an ideal in $\RS_+$
if and only if whenever 
$\alpha \in \CI$ and $\beta \in \RS_+$ so that 
$\alpha \preceq \beta$, then $\beta \in \CI$.

Let $\beta$ be in $\RS_+$. Then $\beta = \sum_{\alpha \in \SI} c_\alpha \alpha$
for $c_\alpha \in \BBZ_{\ge0}$.
The \emph{height} of $\beta$ is defined to be $\hgt(\beta) = \sum_{\alpha \in \SI} c_\alpha$.
Let $\CI \subseteq \RS_+$ be an ideal and let 
$\CI^c := \RS_+ \setminus \CI$
be its complement in 
$\RS_+$. 

\subsection{Arrangement of ideal type}
\label{sect:idealarrangements}
Let $\CA(\RS)$ be the \emph{Weyl arrangement} of $\RS$,
i.e., $\CA(\RS) = \{ H_\alpha \mid \alpha \in \RS_+\}$,
where $H_\alpha$ is the hyperplane in the Euclidean space
$V = \BBR \otimes \BBZ \RS$ orthogonal to the root $\alpha$.
Following \cite[\S 11]{ST06}, 
we associate with an ideal $\CI$ in $\RS_+$ the arrangement 
consisting of all hyperplanes with respect to the roots in $\CI^c$.

\begin{defn}[{\cite[\S 11]{ST06}}]
\label{def:idealtype}
Let $\CI \subseteq \RS_+$ be an ideal.
The \emph{arrangement of ideal type} associated with 
$\CI$ is the subarrangement $\CA_\CI$
of $\CA(\RS)$ defined by 
\[
\CA_\CI := \{ H_\alpha \mid \alpha \in \CI^c\}.
\]
\end{defn}

Thanks to independent fundamental work of Arnold and Saito, the reflection arrangement $\CA(W)$ of any 
real reflection group $W$ is free, cf.~\cite[\S 6]{OT}.
It was shown by Sommers and Tymoczko \cite[Thm.\ 11.1]{ST06} that in case the root system is classical or of type $G_2$, 
each of the subarrangements $\CA_\CI$ of $\CA(W)$ is also free.
The general case was settled only recently in a uniform manner for all types  
in \cite[Thm.\ 1.1]{p-ABCHT-14}.

\begin{theorem}
[{\cite[Thm.\ 11.1]{ST06}, \cite[Thm.\ 1.1]{p-ABCHT-14}}]
\label{thm:ideals}
Let $\RS$ be a reduced root system with 
Weyl arrangement $\CA = \CA(\RS)$.
Let $\CI$ be an ideal in $\RS_+$. 
Then $\CA_\CI$ is free.
\end{theorem}

Moreover, the non-zero exponents of the $\CA_\CI$  
are combinatorially determined, they
are given by the dual of the height partition of the roots in $\CI^c$.

\medskip

It follows from \cite[Thm.\ 5.14]{p-BC10}
that the reflection arrangement $\CA(W)$ of the Weyl group $W$ 
is \emph{inductively free}, see Definition \ref{def:indfree}.
In fact this holds for every Coxeter Group $W$.
The most challenging case here is that of type $E_8$.
Considering the method of proof of 
Theorem \ref{thm:ideals} in \cite{p-ABCHT-14}
and in view of the combinatorial nature of 
the exponents of the free subarrangements $\CA_\CI$ of 
$\CA(W)$, 
it is natural to ask whether the $\CA_\CI$
 are also inductively free.
This question was first investigated in \cite{gR17},
where the following was established in a 
case by case analysis.

\begin{theorem}
[{\cite[Thm.~1.15]{gR17}}]
\label{thm:oldindfree}
Let $\RS$ be a reduced root system with 
Weyl arrangement $\CA = \CA(\RS)$.
Let $\CI$ be an ideal in $\RS_+$. 
Then $\CA_\CI$ is inductively free with the possible 
exception when $W$ is of type $E_8$ 
and $\CI$ is one of $6178$ ideals in $\RS_+$.
\end{theorem}

In essence the proof of Theorem \ref{thm:oldindfree} in \cite{gR17} is based on 
the fact that inductive freeness extends from a localization along a modular 
element in the interesection lattice of corank $1$ to the ambient arrangement, 
see \cite[Thm.~1.12(ii)]{gR17}.
The main result of the present article treats the most difficult case, when $W$ is of type $E_8$,
and thus removes all possible exceptions in Theorem \ref{thm:oldindfree}, confirming
\cite[Conj.~1.16]{gR17}.

\begin{theorem}
\label{thm:indfree}
Let $\RS$ be a reduced root system with 
Weyl arrangement $\CA = \CA(\RS)$.
Let $\CI$ be an ideal in $\RS_+$. 
Then $\CA_\CI$ is inductively free.
\end{theorem}

\begin{remarks}
(a)
Thanks to \cite[Thm.\ 5.14]{p-BC10},
each Weyl arrangement $\CA(\RS)$ is inductively free.
So it is tempting to try to deduce 
inductive freeness for each $\CA_\CI$ by    
using an inductive chain of $\CA(\RS)$ and considering the
induced chain for each $\CA_\CI$.
While this naive approach 
works surprisingly well in some instances
(e.g. for $E_6$),  it does however fail in general. 

(b)
Our proof of the outstanding instances in type $E_8$ is by computational means
and is inspired by methods from \cite[Cor.\ 5.15]{p-BC10}, see \S \ref{s:proof}.
However, the methods from \cite{p-BC10} do not lead to a result for all ideals (in reasonably short runtime):
in \cite{p-BC10} the focus was on exactly one arrangement, $E_8$ (and its restrictions). One of the important steps in \cite{p-BC10} is to guess a ``good'' ordering for the roots. Unfortunately, the ordering proposed in \cite{p-BC10} fails for most of the $25079$ ideals.
As a result, we have to perform approximately 20000 computations of the type treated once in \cite{p-BC10}  (assuming that the small ones are easy).
Since the computation is larger, we need new techniques and computational improvements; see \ref{s:ideas1}, \ref{s:ideas2}, and Remark \ref{remend} for details.\\
Notice further that it is more or less impossible (at the moment) to store \emph{certificates} as was done in \cite{p-BC10}, the amount of data is simply too big (cf.\ Rem.\ \ref{remend}(f)).

(c)
Since inductively free arrangements are 
\emph{divisionally free}, see \cite{tA16}, 
Theorem \ref{thm:indfree} 
affirmatively settles a conjecture by Abe
that all arrangements of ideal type are 
divisionally free, 
\cite[Conj.\ 6.6]{tA16}.

(d)
It would be very desirable to have
a uniform, conceptual
proof of 
Theorem \ref{thm:indfree}.
This would then provide a conceptual 
proof of the fact that the
Weyl arrangement for $E_8$ itself is 
inductively free. 
\end{remarks}

For general information about arrangements, Weyl groups and root systems,  
we refer the reader to \cite{b-BourLie4-6} and 
\cite{OT}.


\bigskip 

\noindent
{\bf Acknowledgments}: 
This work was supported by DFG-grant RO 1072/16-1.

\section{Arrangements and freeness}

\subsection{Hyperplane arrangements}
\label{ssect:arrangements}
Let $V = \BBR^\ell$ be a real $\ell$-dimensional vector space.
A \emph{(real) hyperplane arrangement} $\CA = (\CA, V)$ in $V$ 
is a finite collection of hyperplanes in $V$ each 
containing the origin of $V$.
We also use the term $\ell$-arrangement for $\CA$. 
We denote the empty arrangement in $V$ by $\Phi_\ell$.

The \emph{lattice} $L(\CA)$ of $\CA$ is the set of subspaces of $V$ of
the form $H_1\cap \dotsm \cap H_i$ where $\{ H_1, \ldots, H_i\}$ is a subset
of $\CA$. 
For $X \in L(\CA)$, we have two associated arrangements, 
firstly
$\CA_X :=\{H \in \CA \mid X \subseteq H\} \subseteq \CA$,
the \emph{localization of $\CA$ at $X$}, 
and secondly, 
the \emph{restriction of $\CA$ to $X$}, $(\CA^X,X)$, where 
$\CA^X := \{ X \cap H \mid H \in \CA \setminus \CA_X\}$.
Note that $V$ belongs to $L(\CA)$
as the intersection of the empty 
collection of hyperplanes and $\CA^V = \CA$. 
The lattice $L(\CA)$ is a partially ordered set by reverse inclusion:
$X \le Y$ provided $Y \subseteq X$ for $X,Y \in L(\CA)$.

For $\CA \ne \Phi_\ell$, 
let $H_0 \in \CA$.
Define $\CA' := \CA \setminus\{ H_0\}$,
and $\CA'' := \CA^{H_0} = \{ H_0 \cap H \mid H \in \CA'\}$,
the restriction of $\CA$ to $H_0$.
Then $(\CA, \CA', \CA'')$ is a \emph{triple} of arrangements,
\cite[Def.\ 1.14]{OT}. 

Throughout, we only consider arrangements $\CA$
such that $0 \in H$ for each $H$ in $\CA$.
These are called \emph{central}.
In that case the \emph{center} 
$T(\CA) := \cap_{H \in \CA} H$ of $\CA$ is the unique
maximal element in $L(\CA)$  with respect
to the partial order.
A \emph{rank} function on $L(\CA)$
is given by $r(X) := \codim_V(X)$.
The \emph{rank} of $\CA$ 
is defined as $r(\CA) := r(T(\CA))$.

\subsection{Free hyperplane arrangements}
\label{ssect:free}
Let $S = S(V^*)$ be the symmetric algebra of the dual space $V^*$ of $V$.
Let $\Der(S)$ be the $S$-module of $\BBR$-derivations of $S$.
Since $S$ is graded, 
$\Der(S)$ is a graded $S$-module.

Let $\CA$ be an arrangement in $V$. 
Then for $H \in \CA$ we fix $\alpha_H \in V^*$ with
$H = \ker \alpha_H$.
The \emph{defining polynomial} $Q(\CA)$ of $\CA$ is given by 
$Q(\CA) := \prod_{H \in \CA} \alpha_H \in S$.
The \emph{module of $\CA$-derivations} of $\CA$ is 
defined by 
\[
D(\CA) := \{\theta \in \Der(S) \mid \theta(Q(\CA)) \in Q(\CA) S\} .
\]
We say that $\CA$ is \emph{free} if 
$D(\CA)$ is a free $S$-module, cf.\ \cite[\S 4]{OT}.

If $\CA$ is a free arrangement, then the $S$-module
$D(\CA)$ admits a basis of $\ell$ homogeneous derivations, 
say $\theta_1, \ldots, \theta_\ell$, \cite[Prop.\ 4.18]{OT}.
While the $\theta_i$'s are not unique, their polynomial 
degrees $\pdeg \theta_i$ 
are unique (up to ordering). This multiset is the set of 
\emph{exponents} of the free arrangement $\CA$
and is denoted by $\exp \CA$.

Terao's \emph{addition deletion theorem} plays a 
pivotal role in the study of free arrangements.

\begin{theorem}[{\cite{p-hT-80}, \cite[\S 4]{OT}}]
\label{thm:add-del}
Suppose that $\CA$ is non-empty.
Let  $(\CA, \CA', \CA'')$ be a triple of arrangements. Then any 
two of the following statements imply the third:
\begin{itemize}
\item[(i)] $\CA$ is free with $\exp \CA = \{ b_1, \ldots , b_{\ell -1}, b_\ell\}$;
\item[(ii)] $\CA'$ is free with $\exp \CA' = \{ b_1, \ldots , b_{\ell -1}, b_\ell-1\}$;
\item[(iii)] $\CA''$ is free with $\exp \CA'' = \{ b_1, \ldots , b_{\ell -1}\}$.
\end{itemize}
\end{theorem}

\subsection{Inductively free arrangements}
\label{ssect:indfree}

Theorem \ref{thm:add-del}
motivates the notion of 
\emph{inductively free} arrangements,  see 
\cite{p-hT-80} or 
\cite[Def.\ 4.53]{OT}.

\begin{defn}
\label{def:indfree}
The class $\CIF$ of \emph{inductively free} arrangements 
is the smallest class of arrangements subject to
\begin{itemize}
\item[(i)] $\Phi_\ell$ belongs to $\CIF$, for every $\ell \ge 0$;
\item[(ii)] if there exists a hyperplane $H_0 \in \CA$ such that both
$\CA'$ and $\CA''$ belong to $\CIF$, and $\exp \CA '' \subseteq \exp \CA'$, 
then $\CA$ also belongs to $\CIF$.
\end{itemize}
\end{defn}

\begin{remark}
\label{rem:indfreetable}
It is possible to describe an inductively free arrangement $\CA$ by means of 
an  
\emph{induction table}, cf.~\cite[\S 4.3, p.~119]{OT}.
In this process we start with an inductively free arrangement
and add hyperplanes successively ensuring that 
part (ii) of Definition \ref{def:indfree} is satisfied.
This process is referred to as \emph{induction of hyperplanes}.
This procedure amounts to 
choosing a total order on $\CA$, say 
$\CA = \{H_1, \ldots, H_m\}$, 
so that each of the subarrangements 
$\CA_i := \{H_1, \ldots, H_i\}$
and each of the restrictions $\CA_i^{H_i}$ is inductively free
for $i = 1, \ldots, m$.
In the associated induction table we record in the $i$-th row the information 
of the $i$-th step of this process, by 
listing $\exp \CA_i' = \exp \CA_{i-1}$, 
$H_i$, and $\exp \CA_i'' = \exp \CA_i^{H_i}$, 
for $i = 1, \ldots, m$.
For instance, see 
\cite[Tables 4.1, 4.2]{OT}.

However, note that the notion of ``table'' is misleading: this only encodes the very first layer on the top and is thus mainly useful in dimension three. One also needs such tables for each restriction $\CA_i^{H_i}$ and recursively for the restrictions of the restrictions and so on.
So it is indeed much better to call this data an \emph{induction tree}.

For the same reason, it is very important to note that inductive freeness of a subarrangement of a Weyl arrangement also includes inductive freeness of many restrictions which are not at all subarrangements of Weyl arrangements. Since we know no analogue to ideal subarrangements in these restrictions, these arrangements are much more difficult to handle.
\end{remark}

\section{Restrictions of Weyl arrangements}

\subsection{Crystallographic arrangements}
\label{s:cryst}

We recall some fundamental notions introduced in \cite{p-CH09a} and \cite{p-C10}.
Since we do not need all the details from \emph{loc.~cit.}, 
we concentrate on the essential objects, the root systems\footnote{Notice that the sets of roots presented here are not root systems in the classical sense. In previous work these were also called \emph{root sets} to avoid confusion.}.

\begin{defn}\label{cryarr}
Let $(\CA,V=\mathbb{R}^\ell)$ be a \emph{simplicial arrangement}, i.e.\ the connected components of
$V\setminus \bigcup_{H\in\CA}H$ are open simplicial cones,
and let $R\subseteq V$ be a finite set
such that
$$\CA = \{ \alpha^\perp \mid \alpha \in R\}\quad \text{and} \quad\mathbb{R}\alpha\cap R=\{\pm \alpha\}
\text{ for all } \alpha \in R.$$
For a chamber $K$ of $\CA$ let $\SI^K$ denote
the set of normal vectors in $R$ of the walls of $K$ pointing to the inside.
We call $(\CA,V,R)$ a \emph{crystallographic arrangement} if
\begin{equation}
R \subseteq \sum_{\alpha \in \SI^K} \mathbb{Z} \alpha \quad \text{for all chambers } K.
\end{equation}
\end{defn}

\begin{defn}
A crystallographic arrangement $(\CA,V,R)$ defines \emph{root systems} in a natural way. These are the sets
\[ R^K = \left\{ (a_1,\ldots,a_\ell)\in\mathbb{Z}^\ell \:\: \vline \:\: \sum_{i=1}^\ell a_i \beta_i \in R \right\} \subseteq \mathbb{Z}^\ell \]
where $K$ are the chambers of $\CA$ and $\SI^K=\{\beta_1,\ldots,\beta_\ell\}$.
\end{defn}

\begin{remark}
Of course, these sets $R^K$ depend on the orderings of the elements of $\SI^K$. Fixing this ordering for one chamber, we get canonical orderings for all other chambers, see below.
\end{remark}

Now every root system $R^K$ defines reflections in the following way:

\begin{defn}
Let $\alpha_1,\ldots,\alpha_\ell$ be the standard basis of $\mathbb{Z}^\ell$.
For $1\le i,j\le \ell$, define entries $c^K_{i,j}$ of a matrix $C^K$ by
\[ c^K_{i,j}:=-\max\{ k \mid k\alpha_i+\alpha_j \in R^K \}\quad\mbox{for } i\ne j,
\quad\quad c^K_{i,i}:=2. \]
This matrix is called the \emph{Cartan matrix} of the chamber $K$ and it defines reflections
$\sigma^K_1,\ldots,\sigma^K_\ell$ in $\GL(\mathbb{Z}^\ell)$ via
\[  \sigma^K_i (\alpha_j) = \alpha_j - c^K_{ij} \alpha_i, \quad  i,j=1,\ldots,\ell. \]
\end{defn}

If $K$ and $K'$ are adjacent chambers with $\beta\in R^K\cap R^{K'}$, then $\beta=\alpha_i$ for some $i$ and
$\sigma^K_i(R^K)=R^{K'}$.
For a fixed crystallographic arrangement $(\CA,V,R)$, there may be chambers $K,K'$ with $R^K\ne R^{K'}$ and with different Cartan matrices. But in most cases, the number of different root systems $R^K$ is smaller than the number of chambers.

\begin{example}
(a) If $\CA$ is a Weyl arrangement in $\mathbb{R}^\ell$, $C$ its Cartan matrix, and $R$ its root system, then $(\CA,\mathbb{R}^\ell,R)$ is a crystallographic arrangement. All Cartan matrices $C^K$ are equal to $C$ and all $R^K$ are equal.\\
(b) The arrangements denoted $\CA^k_\ell(2)$ in \cite[\S 6.4]{OT} are crystallographic.
\end{example}

\begin{remark}
A complete classification of crystallographic arrangements was obtained in a series of papers: \cite{p-CH09b}, \cite{p-CH09c}, and finally \cite{p-CH10}.
In rank greater than two, these are the Weyl arrangements, the arrangements $\CA^k_\ell(2)$, 
and $74$ other sporadic arrangements.
\end{remark}

The main reason why crystallographic arrangements are useful in our context is that they form a class of arrangements closed under restrictions (onto elements of their intersection lattice). In particular, every restriction of a Weyl arrangement is crystallographic.
Thus all arrangements that appear in an induction tree of a Weyl arrangement are subarrangements of a crystallographic arrangement.

\subsection{Weyl groupoids}
\label{s:groupoids}

One advantage of knowing that the arrangements we want to consider are crystallographic 
stems from the fact that crystallographic arrangements come with a symmetry structure similar to the Weyl group (and ``equal'' to the Weyl group when it is a Weyl arrangement).

Let $(\CA,V,R)$ be a crystallographic arrangement.
We define a groupoid, i.e.\ a category in which all morphisms are isomorphisms: the objects are the chambers of $\CA$ and the morphisms are compositions of reflections as above. This way, for each pair of (not necessarily adjacent) chambers $K$ and $K'$ we have exactly one morphism (a linear map) $\varphi$ with $\varphi(R^K)=R^{K'}$, and this is the composition of the reflections of a gallery from $K$ to $K'$.

Notice that this definition gives a simply connected groupoid. But if we identify all objects (chambers) with equal root systems, then we possibly get a groupoid with non-trivial automorphism groups at each object; this is the smallest \emph{covering} of the simply connected groupoid.

\begin{example}
If $\CA$ is a Weyl arrangement, then the smallest covering has exactly one object and the automorphism group is the Weyl group.
\end{example}

Let us denote by $\Aut(R^K)$ the automorphism group at object $K$ in the smallest covering.
There are two types of symmetries which are relevant to the computation presented later on:

Assume that $\CA$ is a subarrangement of a crystallographic arrangement $(\hat\CA,\mathbb{R}^\ell,\hat R)$ (for example if it is a subarrangement of a restriction of a Weyl arrangement). Then without loss of generality $\CA=\{\alpha^\perp\mid \alpha\in R\}$ for some $R\subseteq \hat R$, and we may assume that $\hat R=\hat R^K\subseteq\mathbb{Z}^\ell$ for some chamber $K$.
\begin{enumerate}
\item Applying an automorphism $\varphi\in\Aut(\hat R^K)$ to $R$ we obtain a (possibly different) subset $\varphi(R)$ of $\hat R$, but the arrangements $\CA$ and $\{\alpha^\perp\mid \alpha\in \varphi(R)\}$ are isomorphic.
\item Applying a morphism $\varphi$ in the Weyl groupoid to $R$ we obtain a subset $\varphi(R)$ in a (possibly different) root system of $\hat\CA$, but the arrangements $\CA$ and $\{\alpha^\perp\mid \alpha\in \varphi(R)\}$ are isomorphic.
\end{enumerate}

\section{Proof of Theorem \ref{thm:indfree}}
\label{s:proof}

In this section, we explain our algorithmic approach to calculating the inductive trees for the 
arrangements of ideal type $\CA_\CI$ in type $E_8$.
The method is very similar to the one presented in \cite{p-BC10}, but in the end we need some substantial improvements here, mainly because the number of computations is much bigger, but also because the structure is considerably more redundant when 
simultaneously considering all ideals (as opposed to only the Weyl arrangement itself).\\
There are two main tools which we use in the algorithm.

\subsection{Symmetries}
\label{s:ideas1}
The first important idea is to avoid computations which we have already performed.
In the recursion of the computation we constantly produce arrangements for each of which we have to decide whether it is inductively free (and determine its exponents) or not. Since such a computation can take several seconds even for an arrangement of low rank, and since we need to handle millions of arrangements, we need a good method to recognize an arrangement as one that we have already treated. This is where we use the Weyl groupoid:\\
Before we start the computation, we determine the smallest coverings of all Weyl groupoids that will be involved in the calculation (for example of all restrictions of the arrangement of type $E_8$).
Assume now that $\CA$ is an arrangement appearing in the computation, i.e.\ a subarrangement of a crystallographic arrangement $\hat\CA$. Using a morphism of the Weyl groupoid (as in (2) above), we may assume without loss of generality that $\CA$ is a subarrangement in our favorite root system (of $\hat\CA$). In a second step we can decide whether this subarrangement is conjugate to an arrangement we already treated under the automorphism group as in (1) above.

\subsection{Characteristic polynomials}
\label{s:ideas2}
The second idea comes from the following problem and is new compared to the concepts used in \cite{p-BC10}. If at some point we need to treat an arrangement which is not inductively free, then frequently the tree we have to check is quite large; indeed, sometimes it is so vast that a naive computation even for a single case of rank $5$ or $6$ would take a year to complete. Luckily it turns out that often in these cases, the characteristic polynomial of the arrangement does not split into linear factors over $\mathbb{Z}[t]$. Since such an arrangement cannot be free by Terao's Theorem \cite[Thm.~4.137]{OT}, a highly optimized function to compute the characteristic polynomial can, under certain circumstances, avoid a very big amount of computations.

\subsection{The algorithms}
\label{s:computation}

\subsubsection{Initialization}
Before we begin to treat each ideal, we need several preliminaries.

Notice first that we often obtain a subset of a crystallographic arrangement which we only `recognize' once the coordinates of the roots have been permuted.
This is why we need the following.
\algo{CanonicalRootSystem}{$R_+$}{Sort the coordinates of $R_+$ in a canonical way}
{the root system $R_+$ of a crystallographic arrangement}{a root system $R'_+$ which is equal to $R_+$ up to permutation of the coordinates}
{
\item View $R_+$ as a $|R_+|\times \ell$-matrix and sort its columns such that the sum of the coordinates is increasing.
\item If this ordering is not unique, choose the lexicographically smallest such matrix. (This requires us to consider all permutations fixing the sequence of sums.)
\item Return the sorted, permutated copy of $R_+$.
}

To easily compute restrictions, it is useful to know all morphisms in the Weyl groupoid (this avoids a lot of linear algebra and is necessary to identify the root systems):
\algo{MorphismsToSimpleRoots}{$R_+$}{For each $\alpha\in R_+$, compute a morphism mapping to a root system in which $\alpha$ is a simple root (i.e.\ a chamber for which $\alpha$ is a wall)}
{the root system $R_+$ of a crystallographic arrangement}{for each $\alpha$, a linear map}
{
\item Compute all root systems of the crystallographic arrangement.
\item Using the Cartan matrices of each root system, compute morphisms to a fixed canonical one and all automorphism groups.
\item With this information of the Weyl groupoid, we obtain the linear maps as required.
}

We need a function that computes the restriction of a crystallographic arrangement including a map to our preferred canonical root system.
\algo{CrystallographicRestriction}{$R_+,\alpha$}{The `restriction' of $R_+$ to the hyperplane $\alpha^\perp$}
{the root system $R_+$ of a crystallographic arrangement, a root $\alpha\in R_+$}{a map to a canonical root system}
{
\item Compute the restriction of the crystallographic arrangement defined by $R_+$ to $\alpha^\perp$. This is done by applying a morphism to $R_+$ which maps $\alpha$ to a simple root (using {\tt MorphismsToSimpleRoots}). In a second step, we `erase' the coordinate corresponding to $\alpha$ and divide by the greatest common divisors of the vectors.
\item Identify the resulting Weyl groupoid.
\item Compute a map to the canonical chamber and the canonically permuted root system.
}

We can now prepare the whole tree of restrictions for a given root system $R_+$ (for example of type $E_8$).
\algo{AllRestrictions}{$R_+$}{Compute all restrictions of a given Weyl arrangement, including all maps to the canonical objects}
{the root system $R_+$ of a Weyl arrangement}{maps (as sequences of labels) from $R_+$ to each element of the intersection lattice.}
{
\item Iteratively compute restrictions of a given arrangement to each hyperplane:
\item Using {\tt CrystallographicRestriction} and {\tt CanonicalRootSystem} we ensure that we compute everything only up to symmetries, i.e.\ we do not compute the complete intersection lattice of the  Weyl arrangement of $R_+$.
\item We only work with labels, not with root vectors. This way, the main computation (see below) is completely free from linear algebra and only consist of lookups.
}

\begin{example}
If we start with the arrangement of type $E_8$, {\tt AllRestrictions} yields the following canonical objects (up to all symmetries):
\begin{itemize}
\item $1$ restriction of rank $7$ with $91$ hyperplanes.
\item $2$ restrictions of rank $6$ with $63,68$ hyperplanes resp.
\item $3$ restrictions of rank $5$ with $41,46,49$ hyperplanes resp.
\item $6$ restrictions of rank $4$ with $24,25,28,30,32,32$ hyperplanes resp.
\item $8$ restrictions of rank $3$ with $13,13,13,16,17,17,19,19$ hyperplanes resp.
\end{itemize}
Compared with the size of the intersection lattice of type $E_8$, see \cite[Table C.23]{OT}, this is a very small set of cases.

Note that for a reflection group $W$ the possible restricted arrangements are in bijection with conjugacy classes of parabolic subgroups of $W$, so in case of Weyl groups
these are in bijection with subsets of a set of simple roots up to conjugacy.  
Orlik and Terao showed in \cite{orlikterao:free} that each arising restricted arrangement 
is always free in a case-by-case argument.
In  \cite[Cor.~6.1]{douglass:adjoint}, Douglass
gave a uniform proof of this fact.
\end{example}

Of course we also need the set of ideals of our Weyl arrangement (for example of type $E_8$). (Notice that this part is much easier than the other ones, we only reproduce it here for the sake of completeness.)
\algo{Ideals}{$R_+$}{Compute all ideals of a Weyl arrangement}
{the root system $R_+$ of a Weyl arrangement}{all subsets of $R_+$ which are ideals}
{
\item Compute the set $P$ of `principal ideals' and set $J:=P$.
\item Repeat $J \leftarrow J \cup \{ I\cup U \mid U \in P, I \in J\}$ until $J$ remains the same.
\item Return $J$.
}

\subsubsection{The main algorithm}
We now present the main algorithm. Notice that we do not require that the arrangement under consideration is an ideal arrangement; the same algorithm may be used to treat arbitrary subarrangements of crystallographic arrangements (and thus of Weyl arrangements).

\algo{InductiveChain}{$R$,$S$,$e$,$u$,$D$,$N$}{Compute exponents of $R$ if it is inductively free using the Weyl groupoid}
{a set of roots $R$, a subset $S$ of $R$, exponents $e$ of $S$, a label $u$ of a crystallographic arrangement containing $R$, a database $D$ of subarrangements with exponents if inductively free, an invariant $N$ characterizing $R$}{exponents of $R$ if it is inductively free or false otherwise; results of the computation are included into $D$}
{
\item If $|R|=|S|$:
\begin{enumerate}
\item Add $N$ and $e$ to $D$.
\item Add invariants to $D$ given by the orbits of $S$ under the automorphism group. This is used to identify $S$ more quickly if it reappears.
\item Return $e$.
\end{enumerate}
\item If $|S|=0$:
\begin{enumerate}
\item Look for $R$ in $D$ via $N$. If it is included, return the exponents or `false' resp.
\item If the rank of $R$ is in $4,\ldots,7$, then check whether $R$ is conjugate to something known in $D$. Here, we use the invariant given by the orbits of the automorphism group. If it is known, return the exponents or `false' resp.
\item If the rank is at most $5$, compute the characteristic polynomial of $R$. If the roots are not integers, mark this information into $D$ and return `false'.
\end{enumerate}
\item Using the exponents $e$, compute all possible sizes $Y$ of restrictions for the next step in the induction process.
\item For each hyperplane $H$ in $R$ and not in $S$:
\begin{enumerate}
\item Using the precomputed information, map $S$ to the restriction $A:=S^H$.
\item If $|A|\in Y$:
\begin{enumerate}
\item If $A$ has rank $2$ or if $|A|=0$, then let $f$ be the exponents of $A$.
Otherwise, set $f:=${\tt InductiveChain}$(A,\emptyset,\{\{0,\ldots,0\}\},$ label of a crystallographic arrangement containing $A,D,$ invariant of $A)$.
\item If $f$ is not `false' and $f$ is a submultiset of $e$, then:
\begin{enumerate}
\item Compute the exponents $e'$ of $S\cup\{H\}$.
\item Return {\tt InductiveChain}$(R,S\cup\{H\},e',u,D,N)$.
\end{enumerate}
\end{enumerate}
\end{enumerate}
\item If $|S|=0$, then include the information `false' for $(R,N)$ to $D$.
\item Return `false'.
}

To check the inductive freeness for all ideals in, say the Weyl arrangement of type $E_8$, we just call
\[ \text{\tt InductiveChain}(R,\emptyset,\{\{0,\ldots,0\}\},\text{ label of }E_8,D,\text{ invariant of }R) \]
for each ideal subarrangement $R$, starting with an initial empty database $D$.

\begin{remarks}\label{remend}
(a) It is important to keep the same ordering of the roots during the whole computation, in particular also when descending to a restriction. This way we have better chances of encountering a subset we already know.

(b) It is quite expensive to check whether two sets are conjugate under the automorphism group (as required in step 2.\ (2)). To avoid too many of these tests, we compute intersections with orbits and compare these first.
The ranks $4,\ldots,7$ for which we check conjugacy are chosen to be optimal for the special case of $E_8$.

(c) The invariant of a subarrangement is a number: if the subarrangement corresponds to the subset $\{a_0,\ldots,a_n\}$ of labels in the canonical root system, this number is $\sum_i a_i 2^{a_i-1}$.

(d) In contrast to \cite{p-BC10}, it is difficult to store the resulting data in such a way that an independent program can check the results, because of the large number of induction trees. We have implemented the above algorithm within {\sc Magma} \cite{magama} 
since it has a powerful function to check conjugacy (see (b)). But a system like {\sc GAP} would be equally appropriate.

(e) Computing characteristic polynomials turns out to be profitable up to rank $5$. A good algorithm to determine these polynomials could enhance everything:
the problem is that this computation is more or less equivalent to computing the intersection lattice of the arrangement; computing the roots of the polynomials is not an issue.
To accelerate this part, we call an external C-program optimized for arrangements with integral coordinates.

(f) Our main algorithm takes less than an hour of CPU time on a 3.2GHz PC to process all $25079$ ideals in type $E_8$.

(g) It is important that the structure of the database $D$ of subarrangements is in such a way that lookups are easy and fast. We store the data in a four dimensional array: for each rank, for each type of root system, for each size of a subarrangement we have a list of invariants of subarrangements, and for each subarrangement we record whether it is inductively free or not and if so we store its exponents.
\end{remarks}


\bibliographystyle{amsalpha}

\newcommand{\etalchar}[1]{$^{#1}$}
\def\cprime{$'$}
\providecommand{\bysame}{\leavevmode\hbox to3em{\hrulefill}\thinspace}
\providecommand{\MR}{\relax\ifhmode\unskip\space\fi MR }
\providecommand{\MRhref}[2]{%
  \href{http://www.ams.org/mathscinet-getitem?mr=#1}{#2}
}
\providecommand{\href}[2]{#2}

\end{document}